\documentclass{amsart}

\usepackage{Preset}
\usepackage[makeroom]{cancel}
\usepackage{nccmath}

\usepackage{mathrsfs}
\usepackage[latin1]{inputenc}
\usepackage[T1]{fontenc}
\usepackage{tikz}
\usepackage{tikz-cd}
\usepackage[colorlinks=true,linkcolor=blue!50!black,anchorcolor=red,citecolor=blue,filecolor=black,menucolor=black,runcolor=black,urlcolor=black]{hyperref}

\setlist[enumerate,1]{label={\textnormal{(\arabic*)}}}

\numberwithin{equation}{section}

\begin{document}

\title[Wandering continua, renormalization, and Lebesgue measure]{Wandering continua, renormalization, and Lebesgue measure on Julia sets}

\author{Willie Rush Lim}
\address{Department of Mathematics, Brown University, Providence, RI 02912, USA}
\email{willie\_rush\_lim@brown.edu}

\begin{abstract}
Sullivan's celebrated no-wandering-domain theorem asserts that rational maps have no wandering Fatou components. 
In this paper, we establish a measure-theoretic counterpart for the Julia set: every wandering continuum in the Julia set of a finite-type entire function has zero area.
For infinitely renormalizable unicritical polynomials, we also prove that almost every point in the renormalization core (the intersection of the orbits of little Julia sets) is recurrent and contained in the postcritical set. 
The main tool is a conservativity theorem: if $f$ is injective on a forward invariant measurable subset $E$ of the Julia set, then $f|_E$ is conservative and almost every point of $E$ is recurrent and contained in the postsingular set. 
\end{abstract}

\maketitle

\section{Introduction}\label{sec:intro}

Let $f$ be either a degree $d \geq 2$ rational map $\RS \to \RS$ or an entire function $\C \to \C$.
Let $F(f)$ be the Fatou set, the region of stability.
The complement is the Julia set $J(f)$, the set of chaotic points.
Denote by $S(f)$ the set of singular values of $f$ (the closure of the set of critical and asymptotic values), and by
\[
P(f):=\overline{\bigcup_{n\ge0}f^n(S(f))}
\]
the postsingular set of $f$.
We say that $f$ is \emph{of finite type} if $S(f)$ is a finite set, \emph{of infinite type} if otherwise. 
If $f$ is rational, $S(f)$ is the finite set $\CV(f)$ of critical values, and $P(f)$ is the postcritical set.
If $f$ is transcendental entire, we denote the escaping set by 
\[
I(f) = \left\{ z \in \C \ : \ f^n(z) \to \infty \text{ as } n \to \infty \right\}.
\]

Sullivan's no-wandering-domain theorem \cite{Sul85} asserts that for any rational $f$, every connected component of $F(f)$ is pre-periodic.
In his proof, he introduced quasiconformal mappings to one-dimensional complex dynamics and drew a parallel with Ahlfors' finiteness theorem from the theory of Kleinian groups.
This parallel is now known as Sullivan's dictionary, and the use of quasiconformal mappings is now prevalent in the modern developments of complex dynamics.

Recently, Zihao Ye \cite{Ye26a} found an elegant alternative proof of no-wandering-domain theorem that avoids quasiconformal deformations entirely. 
The proof rests on the simple principle that on a hyperbolic Riemann surface with finite area, removing finitely many punctures costs a bounded amount of hyperbolic area. 
The method has very quickly found further applications.
\begin{enumerate}[label=(\roman*)]
\item Ye extended it to transcendental entire and meromorphic functions in various classes \cite{Ye26a, Ye26b}.
\item Prochorov, Rempe and Waterman \cite{PRW26} and, independently, Drach, Pardo-Sim\'on and U\v{c}akar \cite{DPU26} resolved a longstanding open question in transcendental dynamics, showing that transcendental entire functions have no wandering domains with bounded orbits. 
\end{enumerate} 

An immediate consequence of the no-wandering-domain theorem is a complete understanding of the dynamics on the Fatou set.
The Julia set, however, is far more mysterious,
particularly from the measure-theoretic point of view. 
Julia sets of positive area exist \cite{EL87,McM87,BC12}.
The NILF conjecture of McMullen and Sullivan asserts that no rational map (other than a flexible Latt\`es map) carries a measurable invariant line field on its Julia set \cite{McM94, McS98}; it implies the central Density of Hyperbolicity (DoH) conjecture.
For unicritical polynomials $f(z)=z^d+c$, these major conjectures are now reduced to those that are infinitely renormalizable \cite{Hub93}.
The two conjectures are expected to also hold for finite-type entire functions \cite{RvS11}; however, it was long known that invariant line field may exist on $I(f) \cap J(f)$ for infinite-type entire functions $f$ \cite{EL87}.

Motivated by the aforementioned problems, it is natural to ask where exactly the Lebesgue measure lives on the Julia set.
In a recent note \cite{Pet26} by Petersen, it was demonstrated that Ye's hyperbolic area method can be used to show that, for an infinitely renormalizable quadratic polynomial, the intersection of all little filled Julia sets containing the critical value has zero area.
In this paper, we further study the measure-theoretic properties of invariant and wandering subsets of the Julia set.

Throughout this paper, ``Area'' refers to the Lebesgue measure Leb, and ``null'' refers to zero area.
A \emph{continuum} is to a non-empty compact connected metrizable space; it is called \emph{non-degenerate} if it is not a singleton. 
For a continuous dynamical system $f: E \to E$, a subset $X \subset E$ is called \emph{wandering} if $f^m(X)\cap f^n(X)=\emptyset$ for all $m,n \geq 0$ with $m\neq n$.
We will use the notation $\text{orb}(X):=\cup_{n \geq 0} f^n(X)$.

\subsection{On injective forward-invariant subsets}

Our first theorem describes how injectivity of a forward-invariant subset of the Julia set imposes strong measure-theoretic constraints.

\begin{thmx}[Conservativity]
\label{thm:A}
Let $f$ be a rational or entire function. 
Let $E\subset J(f)$ be a measurable set such that
\begin{enumerate}[label=\textnormal{(\alph*)}]
    \item\label{ass:a} $E$ is forward invariant, i.e. $f(E)\subset E$;
    \item\label{ass:b} $f$ is injective on $E$;
    \item\label{ass:c} if $f$ is of infinite type, $E$ is orbitally bounded, i.e. every point in $X$ has a bounded orbit.
\end{enumerate} 
Then,
\begin{enumerate}
\item\label{A1} for every measurable subset $X$ of $E$, the set $X \setminus \cup_{n\geq 1} f^n(X)$ is null; \\ in particular, $E\setminus f(E)$ is null;
\item\label{A2} almost every point in $E$ is recurrent;
\item\label{A3} almost every point in $E$ is in $P(f)$.
\end{enumerate}
\end{thmx}

Item \ref{A1} is equivalent to $(f,E,\textnormal{Leb})$ being a \emph{conservative} system in the sense that every wandering measurable subset of $E$ is null.
It immediately implies:

\begin{corx}[Injective wandering sets]
\label{cor:inj-wandering}
    Let $f$ be a rational or entire function.
    Let $A$ be a wandering measurable subset of $J(f)$ such that $f$ is injective on $f^n(A)$ for all $n
    \geq 0$.
    If $f$ is of infinite type, assume additionally that $A$ is orbitally bounded.
    Then, $A$ is null.
\end{corx}

If $f$ is rational, this follows from Sullivan \cite[Theorem 2]{Sul85}; his quasiconformal proof extends to finite-type entire functions but fails for infinite-type ones because the Teichm\"uller space can be infinite-dimensional.
If $f$ is entire and $A$ has bounded orbit, this corollary follows from \cite[Theorem 1.3]{PRW26}.
In this paper, we avoid using quasiconformal deformation and instead use the hyperbolic area method of Ye \cite{Ye26a}. 

Item \ref{A2} essentially follows from \ref{A1}.
It immediately implies:

\begin{corx}
\label{cor:esc}
    Let $f$ be a finite-type transcendental entire function. 
    Every injective forward-invariant measurable subset of $I(f)\cap J(f)$ is null.
\end{corx}

In other words, Lebesgue measure is not concentrated on the injective part of $I(f)$.
Note that $I(f)$ can have positive area, e.g. sine functions \cite{McM87}; refer to the comprehensive survey by Bergweiler and Rempe \cite{BR25} on properties of $I(f)$.

Item \ref{A3} says that the injective part of $J(f)$ is concentrated in $P(f)$.
It follows from \ref{A2} together with Lemma \ref{lem:attraction} asserting that almost every point in $E$ is attracted to $P(f)$ (together with $\infty$ in the transcendental case); note that the Fatou set is allowed to be empty.
Item \ref{A3} immediately implies:

\begin{corx}
\label{cor:P(f)}
    Let $f$ be a rational or entire function.
    If $P(f) \cap J(f)$ is null, then every measurable subset of $J(f)$ satisfying \ref{ass:a}, \ref{ass:b}, \ref{ass:c} is null. 
\end{corx}

In particular, the corollary applies to geometrically finite maps.
It is expected that every quadratic polynomial has a null postcritical set (c.f. \cite[\S1.5]{Lim26b}), but it may fail for higher degree polynomials in general.

\begin{remark}
    Consider a degree $d \geq 2$ rational map $f$ such that $(f,\RS,\text{Leb})$ is conservative. 
    This occurs if $f$ admits an invariant probability measure in the Lebesgue class, e.g. Latt\`es maps.
    Standard ergodic theory tells us that for any forward invariant measurable subset $E\subset \RS$, $f^{-1}(E) \setminus E$ is null, and if $f$ is injective on $E$, then $E$ must be null.
\end{remark}

\subsection{On wandering continua}

Our second theorem is a measure-theoretic counterpart to the no-wandering-domain theorem for the Julia set.

\begin{thmx}
\label{thm:wandering}
    Every finite-type entire function admits no positive-area wandering continuum in its Julia set.
\end{thmx}

For a polynomial $f$ with locally connected Julia set $J(f)$, every wandering continuum within $J(f)$ is degenerate \cite{BL02}.
For a unicritical polynomial $f(z)=z^d+c$, there are two known examples of the occurence of non-degenerate wandering continua, namely when $f$ admits a neutral periodic point with certain irrational rotation numbers (e.g. non-Brjuno), and when $f$ is infinitely satellite renormalizable with certain unbounded combinatorics.
The latter example will be elaborated in the next subsection.
In both examples, $J(f)$ can have positive area \cite{BC12}, and Theorem \ref{thm:wandering} asserts that the area cannot come from wandering continua.

In the proof of Theorem \ref{thm:wandering}, we show that the topological constraint of being a continuum forces eventual injectivity of $f$ (Proposition \ref{prop:inj}), and then Corollary \ref{cor:inj-wandering} becomes applicable.
This strategy gives us a conditional result for infinite-type entire functions:

\begin{thmx}
\label{thm:wandering-2}
Let $f$ be an infinite-type transcendental entire function and let $K\subset J(f)$ be a wandering continuum.
If the orbit of $K$ is bounded and eventually disjoint from $S(f)$, then $K$ is null.
\end{thmx}

It immediately implies:

\begin{corx}
\label{cor:discrete}
    Let $f$ be a transcendental entire function with discrete $S(f) \cap J(f)$.
    Then, every positive-area wandering continuum in $J(f)$ has unbounded orbit.
\end{corx}

The boundedness assumption in Theorem \ref{thm:wandering-2} cannot be dropped.
Recently, Mart\'{i}-Pete, Rempe, and Waterman \cite{MRW25} proved using an approximation technique that every full continuum $K$ with no interior can be realized as an escaping wandering continuum in the Julia set of an entire function $f$, with $f$ injective on its orbit.
In particular, their construction works for $K$ with positive area.
Furthermore, Evdoridou, Mart\'{i}-Pete, and Rempe \cite{EMPR26} constructed a significant generalization that allows oscillating behavior and any prescribed dynamics on such $K$, not necessarily injective.

\subsection{On infinitely renormalizable maps}

Renormalization is the mechanism behind the universality first observed
in period-doubling cascades. 
A unicritical polynomial $f$ is called renormalizable if some iterate, restricted to a neighborhood of a critical point, is a polynomial-like map in the sense of Douady and Hubbard \cite{DH85}.
This produces small copies of the Julia set inside $J(f)$ and small copies of the Multibrot set inside the parameter plane.
Currently, infinitely renormalizable maps are the frontier of the MLC conjecture on the local connectivity of the Mandelbrot set; MLC implies both NILF and DoH.

Let $f(z)=z^d+c$ be an infinitely renormalizable unicritical polynomial. 
For $n \geq 1$, denote by $J_n(f)$ the little (filled) Julia set of the $n$\textsuperscript{th} renormalization of $f$ containing the critical point $0$. 
Let $q_1<q_2<q_3<\ldots$ be the renormalization periods; for every $n\geq 1$, $f^{q_n}: J_n(f) \to J_n(f)$ is a $d$-to-$1$ covering map branched at $0$ and for $j \in \{1,\ldots, q_n-1\}$, the intersection of $f^j(J_n(f))$ and $J_n(f)$ is at most a single point.
Consider the \emph{renormalization core}
\[
\Rset(f) := \bigcap_{n\geq 1}\ \text{orb}(J_n(f)),
\]
the intersection of the orbit of the little filled Julia sets. 
The core is a compact, disconnected, forward invariant subset of the Julia set.
On it, $f$ is semi-conjugate to the adding machine on $\varprojlim\mathbb{Z}/q_n \mathbb{Z}$.
Moreover, $\Rset(f)$ contains $P(f)$, and every connected component is a wandering continuum intersecting $P(f)$.

\begin{thmx}[The renormalization core]
\label{thm:renorm}
Let $f(z)=z^d+c$ be an infinitely renormalizable unicritical polynomial with some degree $d \geq 2$. 
Then,
\begin{enumerate}
    \item every connected component of $\Rset(f)$ is null,
    \item almost every point in $\Rset(f)$ is recurrent, and
    \item $\Area(\Rset(f))=\Area(P(f))$.
\end{enumerate} 
\end{thmx}

If $f$ is robust \cite{McM94,LP25}, then $\Rset(f)$ is a null Cantor set equal to $P(f)$ and every point in it is recurrent.
In general, every component of $\Rset(f)$ is a wandering continuum that intersects $P(f)$ and it can be non-degenerate; this is the case in the known satellite examples with non-locally connected Julia sets \cite{Mil00, Sor00, Lev11}. 
For the component containing the critical value $c$, item (1) was captured in \cite{Pet26}; as a whole, (1) follows directly from Theorem \ref{thm:wandering}.
Items (2) and (3) are consequences of Theorem \ref{thm:A}.

\subsection{Further discussion}

The results in this paper naturally lead us to the following open problems.

\begin{conj}
\label{conj:renorm}
For every infinitely renormalizable unicritical polynomial $f$, 
\[
\Rset(f) = P(f).
\]
\end{conj}

A positive resolution to this conjecture should follow once geometric \emph{a priori} bounds are obtained for all infinitely renormalizable unicritical polynomials; this is closely related to the MLC conjecture.
Theorem~\ref{thm:renorm} shows at least that the difference $\Rset(f) \setminus P(f)$ is invisible to the Lebesgue measure, just as it is invisible to any invariant probability measure in the satellite case \cite{LP25}.

\begin{ques}
\label{q:D}
If $f$ is a rational or finite-type entire function, does every forward invariant injective measurable subset of $J(f)$ necessarily have zero area?
\end{ques}

By Theorem~\ref{thm:A}, it suffices to consider $E\subset P(f)$.

\begin{ques}
\label{q:E}
Does there exist a rational map $f$ that admits a positive-area wandering continuum within $J(f)$?
\end{ques}

Non-degenerate wandering continua can exist for rational maps even when the Julia set is locally connected \cite{Ro08}.
In the current literature, the only mechanism known to produce positive-area results comes from the Buff-Ch\'eritat Lemma \cite{BC12} stating that certain perturbations
of a Siegel map have the postcritical sets in a small neighborhood of the Siegel disk. 
It is not known whether it can be steered into a single wandering component.
Any hypothetical example would be subject to a number of structural constraints:

\begin{proposition}
\label{prop:rational}
  Let $f$ be a rational map and let $K \subset J(f)$ be a positive-area wandering
  continuum. 
  Then, almost every point in $K$ is attracted to $P(f)$, and for infinitely many $n$, $f$ is non-injective on $f^n(K)$ and $f^{n+1}(K)$ separates $\CV(f)$.
\end{proposition}

The attraction to $P(f)$ follows from Lemma \ref{lem:attraction}; it implies that $f$ must be postcritically infinite.
Non-injectivity follows from \cite[Theorem 2]{Sul85} or Theorem \ref{thm:A}.
The separation property follows from Lemma \ref{lem:inj} and Corollary \ref{cor:inj-wandering}.

\subsection{Acknowledgements}
I would like to thank Jeremy Kahn and Lasse Rempe for helpful comments.
This paper is inspired by the new development on wandering domains spurred by Zihao Ye in collaboration with generative artificial intelligence (AI). 
Just like \cite{Pet26}, this paper is produced and written without the use of generative AI.

\section{The proof}

Throughout this section, $f$ is either a degree $\geq 2$ rational map or a non-linear entire function.

\subsection{Conservativity}
\label{ss:prelim}

The following is the key lemma to prove Theorem \ref{thm:A}.

\begin{lemma}[{\cite{Ye26a}, \cite{PRW26}, \cite[Lemma~1]{Pet26}}]\label{lem:gap}
Let $E \subset \RS$ be a measurable set satisfying \ref{ass:a}, \ref{ass:b}, and \ref{ass:c}. 
If $E \setminus f(E)$ is contained in $J(f)$, then $E \setminus f(E)$ is null.
\end{lemma}

Below, for any domain $U \subset \RS$ that misses at least three points on $\RS$, denote by $\lambda_U = \lambda_U(z)|dz|$ the complete hyperbolic metric on $U$.

\begin{proof}
For rational $f$ and for $n \geq 1$, let $X_n$ be the set of repelling periodic points of period at most $n$ and add the critical values to get $Y_n:=X_n\cup\CV(f)$.
Also, denote $\rho_n:=\lambda_{\RS \setminus X_n}$ and $\sigma_n:=\lambda_{\RS \setminus Y_n}$. 
The argument of \cite{Ye26a} (see also \cite[Lemma~1]{Pet26}) gives
\[
\int_E\rho_n^2\,dA\le\int_{f(E)}\sigma_n^2\,dA\le2\pi(2d-2)+\int_{f(E)}\rho_n^2\,dA.
\]
The first inequality uses Schwarz-Pick Lemma applied to the covering map $f: \RS \setminus f^{-1}(Y_n) \to \RS \setminus Y_n$, and the injectivity of $f$ on $E$;
the second uses Gauss-Bonnet and $|\CV(f)|\leq 2d-2$. 
By forward invariance, 
\[
\int_{E\setminus f(E)}\rho_n^2\,dA \leq 2\pi(2d-2) \qquad \text{ for all }n \geq 1.
\]
The density $\rho_n$ increases to $\infty$ at every point of $J(f) \setminus \cup_n X_n$.
By monotone convergence, $E\backslash f(E)$ must be a null set.

If $f$ is entire of finite type, the same argument applies once $\infty$ is added to the punctures; compare \cite{Ye26a}, where this is done for wandering domains. 
Here $X_n$ is again the set of repelling periodic points of period at most $n$, which are dense in $J(f)$, and $Y_n:=X_n\cup S(f)$. 
Then, take $\rho_n:=\lambda_{\C\setminus X_n}$ and $\sigma_n:=\lambda_{\C\setminus Y_n}$.
The same chain of inequalities apply, with the bound $2\pi(2d-2)$ replaced by $2\pi|S(f)|$.

Finally, let $f$ be transcendental entire of infinite type. 
For $R>0$, let
\[
E_R:=\{z\in E:\ |f^n(z)| \leq R\ \text{for all } n \geq 0\}.
\]
Each $E_R$ is measurable, bounded, forward invariant, and $f$ is injective on $E_R$. 
Assuming orbital boundedness, the union of $E_R$ over all $R \in \N$ is equal to $E$. 
It therefore suffices to show that the set $D_R := E_R \backslash f(E)$ is null for every $R$.
From now on, let us fix $R$.

Let $V:=\{|z|<2R\}$, and let $L :=f(\CP(f)\cap V)$, which is a finite set. 
By \cite[Lemma~4.2]{PRW26}, there exist finite sets $P_1\subset P_2\subset\cdots\subset\overline V$ with $f(P_j\cap V)\subset P_j$. 
These $P_j$'s are obtained by pulling back, within $V$, finite subsets of $\partial V$ whose union is dense in $\partial V$. 
The closure of $\bigcup_j P_j$ contains $T_V(f)\setminus\interior T_V(f)$, where $T_V(f)$ is the set of points whose orbits stay in $V$.
By removing from $E_R$ the null set of points whose orbits meet $\cup_{j \geq 1} f^{-1}(P_j\cup L)$, we may assume that $E_R$ avoids $f^{-1}(P_j\cup L)$ for every $j$.

Since $f$ is not a covering, Schwarz--Pick is replaced by the estimate \cite[Theorem~3.1]{PRW26} on the compact disk $\{|z|\leq R\}$. 
It bounds the defect between $\lambda_{\C\setminus P_j}$ and the pullback of $\lambda_{\C\setminus(P_j\cup L)}$ by a constant $C=C(R)>0$ independent of $j$. 
The resulting inequalities \cite[Lemma~4.1]{PRW26} give
\[
\Area_{\C\setminus P_j}\bigl(D_R\bigr)
\leq 
\Area_{\C\setminus P_j}\bigl(E_R\setminus f(E_R)\bigr)
\leq 
C+2\pi|L|
\]
for all $j \geq 1$.
Points of $D_R$ lie in $J(f)$ and have orbits in $V$. 
By Montel's theorem, the interior of $T_V(f)$ is in $F(f)$, so $D_R$ is contained in $T_V(f)\setminus\interior T_V(f)$. 
Hence, $\lambda_{\C\setminus P_j}\to\infty$ on $D_R$ outside the countable set $\cup_j P_j$ \cite[Lemma~2.4]{PRW26}. 
Therefore, $D_R$ is a null set.
\end{proof}

Next, let $\hat{P}(f) := P(f) \cup \{\infty\}$ if entire, $\hat{P}(f) := P(f)$ if rational.
The following lemma is inspired by McMullen's Ergodic vs. Attracting Theorem \cite{McM94}.

\begin{lemma}[Attraction to $\hat{P}(f)$]
\label{lem:attraction}
Let $E \subset J(f)$ be a forward invariant measurable set such that $\RS \backslash E$ has positive area. 
Then, for almost every $z\in E$,
\[
\dist \bigl( f^n(z), \hat{P}(f) \bigr) \to 0 \qquad \text{as } n \to \infty.
\]
\end{lemma}

\begin{proof}
Assume that $\hat{P}(f)$ is not the whole sphere; otherwise, there is nothing to prove.
By conjugating with a M\"obius transformation, assume that $\infty \in \hat{P}(f)$.
If $\widehat P(f) \cap \C$ were a singleton, say $\{s\}$, then $f$ would be a covering map from $\C\setminus f^{-1}(s)$ onto the punctured plane $\C\setminus\{s\}$.
This would imply that $f$ is rational and $f$ is affinely conjugate to $z\mapsto z^{\pm d}$. 
In that case, $J(f)$ is a circle and there is nothing to prove. 
Hence, we will assume that $\hat{P}(f) \cap \C$ contains at least two points.

Suppose that there exists a positive-area forward invariant subset $E \subset J(f)$ such that $\dist(f^n(z),\hat{P}(f)) \not \to 0$ for every $z\in E$. 
Below, we will show that $E$ has full Lebesgue measure.

Let $z$ be a Lebesgue density point of $E$. 
There exist times $n_i \to \infty$, a point $y \in \C \setminus P(f)$, and $\delta>0$ such that $\D(y,3\delta) \cap P(f) = \emptyset$ and that the point $w_i:=f^{n_i}(z)$ lies in $\D(y,\delta)$ and converges to $y$.
In particular, $\D(w_i,2\delta) \subset \D(y,3\delta)$.

Since $\D(w_i,2\delta)$ is disjoint from $P(f)$, there is a univalent branch $g_i$ of $f^{-n_i}$ on $\D(w_i,2\delta)$ with $g_i(w_i)=z$.
By Montel, the maps $g_i$ form a normal family because they omit $\hat{P}(f)$. 
Any subsequential limit of $g_i$ is constant; otherwise, $f^{n_i}$ would be normal on an open set containing $z\in J(f)$. 
Hence, $U_i:=g_i(\D(w_i,\delta))$ are shrinking disk neighborhoods of $z$.

Set $B:=\D(y,\delta/2)$, an open disk neighborhood of $y$.
We have
\begin{align*}
    \frac{\Area(B\backslash E)}{\Area \D(w_i,\delta)} &\leq \frac{\Area\bigl(\D(w_i,\delta)\setminus f^{n_i}(E\cap U_i)\bigr)}{\Area \D(w_i,\delta)} =
\frac{\Area\bigl(f^{n_i}(U_i\setminus E)\bigr)}{\Area\bigl(f^{n_i}(U_i)\bigr)} \\
&\preceq \,
\frac{\Area(U_i\setminus E)}{\Area(U_i)} 
\longrightarrow 0.
\end{align*}
The first inequality is because $B\subset \D(w_i,\delta)$ and $f^{n_i}(E\cap U_i)\subset E$,
the third follows from Koebe distortion theorem for $f^{n_i}|_{U_i}$, and the convergence to $0$ is because $z$ is a density point of $E$.
Hence, $\Area(B\setminus E)=0$.

The set $B\setminus J(f)$ is open and has zero area because $E\subset J(f)$. 
Hence it is empty and $B\subset J(f)$. 
By Montel, the orbit of $B$ has full measure. 
For every $n \geq 1$, the set $f^n(B)\setminus E$ is null because it is contained in the null set $f^n(B\setminus E)$. 
Hence, $E$ has full measure. 
\end{proof}

\begin{lemma}
\label{lem:attraction-2}
    The map $f$ cannot be injective on any set of full Lebesgue measure.
\end{lemma}

\begin{proof}
    Let $E$ be any set of full measure. 
    Let us fix non-critical points $a\neq b$ with $f(a)=f(b)$.
    Choose disjoint open disk neighborhoods $A\ni a$ and $B\ni b$ such that $f$ maps both $A$ and $B$ univalently onto the same open disk $D=f(A)=f(B)$.
    Observe that $E \cap A$ has full measure in $A$ and $E \cap B$ has full measure in $B$. 
    Then, $f(E \cap A)$ and $f(E \cap B)$ both have full measure in $D$, so they intersect. 
    A point in the intersection has preimages in both $A$ and $B$, which lie in $E$ and are distinct. 
    Hence, $f$ is not injective on $E$.
\end{proof}

We are now ready to prove Theorem \ref{thm:A}.

\begin{proof}[Proof of Theorem \ref{thm:A}]
Let $f$ be rational or transcendental entire, and let $E$ be a measurable subset of $J(f)$ that is forward invariant \ref{ass:a}, injective \ref{ass:b}, and bounded \ref{ass:c} if of infinite type.

Let $X$ be a measurable subset of $E$.
Then, $\text{orb}(X)$ is a measurable subset of $J(f)$ satisfying \ref{ass:a}, \ref{ass:b}, and \ref{ass:c}.
By Lemma~\ref{lem:gap}, $\text{orb}(X) \setminus f(\text{orb}(X)) = X \setminus \cup_{n \geq 1}f^n(X)$ is null.

Next, let $\{B_i\}_{i\ge1}$ be a countable basis of open sets in $\RS$ and let
\[
Y_i:= \{z\in E\cap B_i\;:\;f^n(z)\notin B_i\text{ for all }n\ge1\}
\qquad
\text{for } i \geq 1.
\]
Observe from the definition that for every $n \geq 1$, $Y_i\cap f^n(Y_i)$ always empty. 
By the previous paragraph, $Y_i$ is a null set. 
Every point $z$ in $E \setminus \cup_i Y_i$ returns to every basic open set containing it.
Hence, almost every point in $E$ is recurrent.

By the injectivity on $E$ and by Lemma \ref{lem:attraction-2}, the complement of $E$ has positive area.
Then, by Lemma \ref{lem:attraction}, a typical point $z$ in $E$ converges to $\hat{P}(f)$ under iteration. 
But since $z$ is also typically recurrent, then $z$ is contained in $\hat{P}(f)$.
In the transcendental case, $z\neq\infty$ and so $z$ is also in $P(f)$.
Thus, $E\setminus P(f)$ is null.
\end{proof}


\subsection{Wandering continua}
\label{ss:continua}

In this subsection, we will prove Theorems \ref{thm:wandering} and \ref{thm:wandering-2}.
Recall that a continuum $K \subset \RS$ is called \emph{full} if its complement is connected, and the \emph{fill} of a continuum $K \subset \C$ is the smallest full continuum in $\C$ containing $K$.

The two lemmas below use the no-wandering domain theorem.

\begin{lemma}
\label{lem:full}
Let $K\subset J(f)$ be a wandering continuum with bounded orbit. 
Then, for all $m \in \N$, $f^m(K)$ is full.
\end{lemma}

\begin{proof}
Suppose for a contradiction that the complement of $f^m(K)$ has a bounded component $U$ for some $m \in \N$. 
For every $j\ge0$, $f^j(U)$ is a bounded connected open set with boundary contained in $\overline{\text{orb}(K)}$. 
Thus, each $V_j$ is contained in the filled hull of $\overline{\text{orb}(K)}$.
By Montel, $U$ has to lie in the Fatou set. 
Since $\partial U\subset J(f)$, then $U$ must be a component of $F(f)$.
Since the orbit of $U$ is bounded, the no-wandering-domain theorem \cite{Sul85, Ye26a, PRW26, DPU26} implies that $U$ must be pre-periodic with some pre-period $p\geq 0$ and period $q\geq 1$.
Then, the non-empty set $\partial V_{p}$ lies in both $f^{m+p}(K)$ and $f^{m+p+q}(K)$, contradicting the wandering property.
\end{proof}

\begin{lemma}
\label{lem:bdd-comp}
    Suppose $f$ is a finite-type entire function and $K\subset J(f)$ is a wandering continuum.
    Any bounded connected component of $\C \backslash K$ has to intersect $J(f)$.
\end{lemma}

\begin{proof}
    Suppose for a contradiction that there is a bounded connected component $U$ of $\C \backslash K$ that is contained in $F(f)$.
    Then, $U$ is Fatou component because $\partial U \subset K \subset J(f)$.
    By the no-wandering-domain theorem \cite{Sul85,Ye26a}, $U$ is pre-periodic.
    Just like the previous lemma, this would contradict the wandering property for $K$.
\end{proof}

The next two lemmas allow us to deduce injectivity out of topological constraints.

\begin{lemma}
\label{lem:inj}
Let $K' \subset \C$ be a connected set with image contained in a full continuum $K$ disjoint from $S(f)$. 
Then, $f$ is injective on $K'$.
\end{lemma}

\begin{proof}
Assume without loss of generality that $K$ is non-degenerate disjoint from $\infty$.
The hypothesis implies the existence of a small open disk neighborhood $W$ of $K$ disjoint from $S(f)$.
(For example, one can systematically construct $W$ using the Riemann mapping of the complement of $K$.)
Since $f:f^{-1}(W)\to W$ is a covering map and $W$ is simply connected, every component of $f^{-1}(W)$ is mapped homeomorphically onto $W$. 
The connected set $K'$ lies in a single component of $f^{-1}(W)$, and we are done.
\end{proof}

\begin{lemma}
\label{lem:L-set}
    Let $f$ be transcendental entire and let $L \subset \C$ be a closed connected set containing $S(f) \cup \{0,f(0)\}$.
    Let $K \subset \C$ be a continuum avoiding $L \cup f^{-1}(L)$ and $K' \subset \C$ be a connected set such that $f(K') \subset K$.
    Then, $f$ is injective on $K'$.
\end{lemma}

In application, $L$ will be taken to be a closed disk. 
This lemma is inspired by the construction of logarithmic tracts for bounded-type entire functions; see the survey by Sixsmith \cite{Six18}.

\begin{proof}
    We first claim that every component of $\C \setminus f^{-1}(L)$ is simply connected.
    Let $\Omega$ be a component of the connected closed set $L$.
    There are two cases: \vspace{0.03in}
    
    \noindent\underline{Case 1:} Either $L$ is unbounded or $\Omega$ is bounded. \\
        Then, $\Omega$ is simply connected and $f$ sends each component of $f^{-1}(\Omega)$ onto $\Omega$ as a conformal isomorphism. \vspace{0.03in}
        
    \noindent\underline{Case 2:} $L$ is bounded and $\Omega$ is unbounded. \\
        Then, $\Omega$ is doubly connected.
        Since $f$ is transcendental, each component of $f^{-1}(\Omega)$ maps onto $\Omega$ as an infinite degree universal cover. \vspace{0.03in}
    
    Since $K$ is a connected set avoiding $f^{-1}(L)$, it is contained in a unique component $T$ of $\C \setminus f^{-1}(L)$.
    Since $T$ is simply connected, the fill $\tilde{K}$ of $K$ is contained in $T$, hence disjoint from $f^{-1}(L)$.
    The assumption that $L$ contains $\{0,f(0)\}$ implies that $L$ intersects $f^{-1}(L)$.
    Together with the assumption that $K$ avoids $L$, then the fill $\tilde{K}$ is also disjoint from $L$ and in particular $S(f)$.
    Then, by Lemma \ref{lem:inj}, $f$ is injective on $K'$.
\end{proof}

Theorems \ref{thm:wandering} and \ref{thm:wandering-2} follow directly from Corollary \ref{cor:inj-wandering} and the following proposition.

\begin{proposition}[Injectivity of wandering continua]
\label{prop:inj}
    Let $f$ be an entire function and $K \subset J(f)$ be a wandering continuum.
    If $f$ is of infinite type, assume that the orbit of $K$ is bounded and eventually disjoint from $S(f)$.
    Then, $f$ is eventually injective along the orbit of $K$.
\end{proposition}

\begin{proof}
    Suppose $f$ is of infinite type.
    By the additional hypothesis, Lemmas~\ref{lem:full} and \ref{lem:inj} directly imply that $f$ is eventually injective along the orbit of $K$.
    Thus, we will suppose from now on that $f$ is of finite type.
    Assume without loss of generality that each $K_n := f^n(K)$ is disjoint from $S(f)$.
    We will suppose for a contradiction that there exists infinitely many times $n_1 < n_2 < n_3 < \ldots $ such that $f^{n_j} : K_{n_j} \to K_{n_{j+1}}$ is not injective.
    
    For every $n \geq 0$, let $\tilde{K}_n$ be the fill of $K_n$.
    By the maximum principle, $f(\tilde{K}_n) \subset \tilde{K}_{n+1}$.
    By Lemma \ref{lem:inj} and the non-injectivity of $f$ on each $K_{n_j}$, the fill $\tilde{K}_{n_i+1}$ would have to contain a singular value.
    But since $K_{n_j+1}$ does not contain any singular value, it would mean that the complement of $K_{n_j+1}$ has a bounded connected component $V_{j}$.

    By Lemma \ref{lem:bdd-comp}, for $j \in \N$, each $V_j$ meets $J(f)$.
    Let $\mathbf{R} >0$ be such that $\{|z| \geq \mathbf{R}\}$ is free of any exceptional value.
    The blow-up property of $J(f)$ implies that for any $R\geq \mathbf{R}$, there exists some $N_R \in \N$ such that for $k \geq N_R$, $f^k(V_j)$ contains the circle $\{|z| = R\}$; consequently, the fill $\tilde{K}_{n_j+k}$ contains the whole closed disk $\{|z| \leq R\}$.
    Let $R_j := \max \bigl\{j, \max_{z \in K_{n_j+1}} |z|, \mathbf{R} \bigr\}$.
    By passing to further subsequence, assume that $n_{j+1} > n_j + N_{R_j}$.
    Then, for $k > n_{j+1}$, $\tilde{K}_k$ contains $\tilde{K}_{n_j+1}$ which contains $\{|z| \leq j-1\}$.
    Since $K_k$ and $K_{n_j+1}$ are disjoint, we conclude that for all $j \geq 1$ and $k> n_{j+1}$, the continuum $K_k$ lies in $\{|z|>j-1\}$.
    
    What we have just shown above is that the orbit of $K$ has to escape to infinity uniformly under $f$.
    If $f$ is a polynomial, this would violate the boundedness of $J(f)$.
    If $f$ is transcendental, then the orbit of $K$ must eventually avoid a closed disk containing $S(f) \cup \{0,f(0)\}$, and by Lemma \ref{lem:L-set}, $f$ would have to be eventually injective along the orbit of $K$.
\end{proof}


\subsection{The renormalization core}
\label{ss:renorm}

Let $f(z)=z^d+c$ be an infinitely renormalizable unicritical polynomial of some degree $d \geq 2$.
It admits an infinite sequence of renormalizations in the sense of Douady-Hubbard polynomial-like maps \cite{DH85} with growing periods $q_1 < q_2 < q_3 < \ldots$ with $q_{n+1}$ divisible by $q_n$. 
We will set $q_0=1$ for convenience.
For each $n \geq 1$, let $J_n = J_n(0)$ be the little (filled) Julia set about the critical point $0$ for the $n$\textsuperscript{th} renormalization, and let $J_n(j) = f^j(J_n)$ for $j \geq 1$.
Denote
\[
\mathcal{J}_n:= \bigcup_{j=0}^{q_n-1} J_n(j), \qquad 
\Rset=\bigcap_n \mathcal{J}_n, \qquad 
J_\infty:=\bigcap_n J_n.
\]
Before we prove Theorem \ref{thm:renorm}, we list a number of basic facts about these compact sets.
A thorough exposition can be found in the book of McMullen \cite[Ch.~7]{McM94} in the quadratic case; the higher degree case is analogous, see e.g. \cite{Ino02}.
\begin{enumerate}[label=\textnormal{(\roman*)}]
    \item The filled Julia set of $f$ has empty interior and is thus equal to $J(f)$.
    \item The map $f$ gives us the following sequence of branched covering maps:
\[
J_n(0) \ \xrightarrow[\ \ d:1 \ \ ]{} \ J_n(1) \ \xrightarrow[\ \ 1:1 \ \ ]{} \ J_n(2) \ \xrightarrow[\ \ 1:1 \ \ ]{} \  \ldots \ \xrightarrow[\ \ 1:1 \ \ ]{} \  J_n(q_n) = J_n(0).
\]
    \item $J_{n+1}(i)$ is contained in $J_n(j)$ if and only if $i\equiv j$ mod $q_n$. 
    \item If any two distinct little Julia sets $J_n(i)$ and $J_n(j)$ intersect, then $i \equiv j$ mod $q_{n-1}$ and the intersection is a repelling periodic point of period $q_{n-1}$.
\end{enumerate}
These standard facts imply the following properties.
\begin{enumerate}[label=\textnormal{(\Roman*)}]
    \item $J_\infty$ is the connected component of $\Rset$ containing $0$.
    \item $f: \Rset \backslash J_\infty \to \Rset \backslash f(J_\infty)$ is a homeomorphism. 
    \item Every connected component of $\Rset$ is a wandering continuum.
\end{enumerate}
Then, Theorem~\ref{thm:renorm} becomes an immediate application of Theorems \ref{thm:A} and \ref{thm:wandering}.

\begin{proof}[Proof of Theorem~\ref{thm:renorm}]
By Theorem~\ref{thm:wandering} and (III), every connected component of $\Rset$ is a null set. 
By (I), the set $E = \Rset \backslash (\cup_{j\geq 0} f^{-j} (J_\infty))$ is a null set.
By (II), $E$ is forward invariant and $f$ is injective on $E$.
Then, apply Theorem~\ref{thm:A} to $E$.
We obtain that almost every point in $\Rset$ is recurrent and $\Area(\Rset)=\Area(P(f))$.
\end{proof}

\bibliographystyle{alpha}
 
{\small \bibliography{bibliography}}

\end{document}